\documentclass[11pt]{amsart}
\usepackage{geometry}                
\usepackage{graphicx}
\usepackage{amssymb}
\usepackage{amscd}
\usepackage{epstopdf}
\font\tency=wncyr10  \def\sf{\hbox{\tency X}} 

\newtheorem{df}{Definition}[section]
\newtheorem{conjecture}{Conjecture}[section]
\newtheorem{thm}{Theorem}[section]
\newtheorem{prop}{Proposition}[section]
\newtheorem{lm}{Lemma}[section]

\newtheorem{fact}{Fact}[section]

\title{An archimedian analog of Iwasawa theory}
\author{Ken-ichi Sugiyama,}
\thanks{Research supported by in part by JSPS grants Kiban(C)22540068\\
Address : Chiba University, Yayoi-cho Inage-ku, Chiba, Japan \\
E-mail address : sugiyama@math.s.chiba-u.ac.jp
}

\begin{document}
\maketitle
\begin{abstract}
We will show a conjecture which reduces Mazur-Tate-Teitelbaum conjecture to the known cases. In order to explain its background  we will develop an archimedian analog of Iwasawa theory. Moreover consequences of the conjecture which are related to Birch and Swinnerton-Dyer conjecture will be discussed.
AMS classification 2000: 11F11, 11F67, 11F85, 11G05, 11G40
 \end{abstract}
\section{Introduction}
In this report we will describe a conjecture which reduces Mazur-Tate-Teitelbaum conjecture (see {\bf Conjecture 3.1}) to the known cases.
We fix a prime  $p$ greater than or equal to $5$. Let $E$ and $E^{\prime}$ be elliptic curves defined over ${\mathbb Q}$. We say that they have ordinary reduction {\rm of the same type} at $p$ if one of the following conditions holds: 
\begin{enumerate}
\item they have good ordinary reduction  at $p$ and the cardinality of ${\mathbb F}_p$-rational points of their reductions are equal,
\item they have split multiplicative reduction at $p$,
\item they have non-split multiplicative reduction at $p$. 
\end{enumerate}
In particular if $E$ has ordinary reduction at $p$ and if their reduction $\tilde{E}_p$ and $\tilde{E}^{\prime}_p$ are isomorphic over ${\mathbb F}_p$ they have the same type. 
\begin{conjecture}
Let $E$ and $E^{\prime}$ be elliptic curves defined over ${\mathbb Q}$. Suppose that they have ordinary reduction of the same type at $p$. Then
\[{\rm ord}_{s=1}L(E,s)-{\rm ord}_{s=0}{\mathcal L}_{E,p}(s)={\rm ord}_{s=1}L(E^{\prime}, s)-{\rm ord}_{s=0}{\mathcal L}_{E^{\prime},p}(s).\]
\end{conjecture}
It is clear that Mazur-Tate-Teitelbaum conjecture implies {\bf Conjecture 1.1}. Conversely we will show that {\bf Conjecture 1.1} indues theirs. Here is a brief outline of an argument (see \$3.1 for details). It is known that Mazur-Tate-Teitelbaum conjecture is true for an elliptic curve $E$ defined over ${\mathbb Q}$ whose L-function does not vanish at $s=1$. More precisely if $E$ has a good ordinary or a non-split multiplicative reduction at $p$ it is obviously true by the definition (see {\bf Fact  2.1}). If $E$ has a split multiplicative reduction it is a theorem due to Greenberg-Stevens \cite{Greenberg-Stevens} (see also \cite{Kobayashi2006}). Suppose that our conjecture were true. It is sufficient to find an elliptic curve $E^{\prime}$ defined over ${\mathbb Q}$ which has an ordinary reduction of the same type at $p$ as $E$ and that $L(E^{\prime},1)\neq 0$. Using results of Ono and Skinner \cite{Ono-Skinner1998} we will construct such a curve by a quadratic twist of $E$.\\
 
Let us briefly explain a motivation of the conjecture. Let  $\Gamma_{\infty}$ be the set of $p$-adic integers congruent one modulo $p$. For an elliptic curve defined over ${\mathbb Q}$ which has ordinary reduction at $p$, Mazur, Tate and Teitelbaum constructed an element $\mu_{E,p}\in {\mathbb Z}_{p}[[\Gamma_{\infty}]]$  which interpolates special values of twisted L-function of $E$ at $s=1$ (see {\bf Fact 2.1}). The $p$-adic L-function of $E$ is intuitively 
\[{\mathcal L}_{E,p}=\chi_{s}(\mu_{E,p}),\]
where $\chi_{s}$ is a character of $\Gamma_{\infty}$ defined by $\chi_{s}(x)=x^{-s}$.
Let $E$ and $E^{\prime}$ be elliptic curves defined over ${\mathbb Q}$ which have ordinary reduction of the same type at $p$ and $\chi$ a character of $\Gamma_{\infty}$ of finite order. Then
\[\sigma(\chi(\phi(\mu_{E,p})))\frac{L(E^{\prime},\chi^{1-p},1)}{\Omega_{E^{\prime}}}=\sigma(\chi(\phi(\mu_{E^{\prime},p})))\frac{L(E,\chi^{1-p},1)}{\Omega_{E}},\]
where $\phi$ is a homomorphism of ${\mathbb Z}_{p}[[\Gamma_{\infty}]]$ induced by an automorphism $x \mapsto x^{p-1} \,(x \in \Gamma_{\infty})$ and $\sigma$ is an isomorphism from ${\mathbb C}_{p}$ to ${\mathbb C}$. 
Suppose we were able to constructed a ${\mathbb C}$-valued measure $\xi_{\infty,E}$ and $\xi_{\infty,E^{\prime}}$ on $\Gamma_{\infty}$ satisfying
\begin{enumerate}
\item $\chi(\xi_{\infty,E})=L(E,\chi^{1-p},1)$ for a finite character of $\Gamma_{\infty}$,
\item $\chi_{s}(\xi_{\infty,E})=L(E, 1+(1-p)s)$,
\end{enumerate}
and so does $\xi_{\infty,E^{\prime}}$.
If we brutely replace $\chi$  in the above equation by $\chi_{s}$ we will obtain
\[\sigma({\mathcal L}_{E,p}((p-1)s))\frac{L(E^{\prime},1+(1-p)s)}{\Omega_{E^{\prime}}}=
\sigma({\mathcal L}_{E^{\prime},p}((p-1)s))\frac{L(E,1+(1-p)s)}{\Omega_{E}},\]
which implies {\bf Conjecture 1.1}.  In order to realize this idea we will develop a ${\mathbb C}$-valued measure theory on $\Gamma_{\infty}$ for an elliptic curve in \$2.2 and \$2.3, which is an arichmedian analog of Iwasawa theory. The motivation will be explained in \$2.4.  (Note that the naively  conjectured this equation is seemed to be too strong. In fact if $E$ and $E^{\prime}$ has a split multiplicative reduction at $p$ and if their L-functions does not vanish at $s=1$ it says that their L-invariants should be equal. But this fact is true if $E^{\prime}$ is a quadratic twist of $E$, which we may impose to derive consequences from {\bf Conjecture 1.1}.)  In \$3.1 we will discuss a relation between {\bf Conjecture 1.1} and the conjecture of Mazur-Tate-Teitelbaum. The remaining sections will be devoted to an application to Birch and Swinnerton-Dyer conjecture.\\



Throughout the paper we will use the following notation. $p$ will be a prime greater than or equal to $5$. Let $G$ be a group and $K$ be a commutative field. We denote the group ring of $G$ whose coefficients are in $K$ by $K[G]$. Fixing an embedding we will consider $\bar{\mathbb Q}$ as a subfield of ${\mathbb C}$ and ${\mathbb C}_p$.
\section{A motivation of the conjecture}
\subsection{Review of the theory of $p$-adic measures}
In this section we will review the theory of $p$-adic integrals and $p$-adic L-functions. Details will be found in \cite{Hida}, \cite{MS} and \cite{MTT}. 
For a positive integer $r$ let $\Gamma_{r}$ be the kernel of the mod $p$ reduction map,
\[({\mathbb Z}/(p^{r+1}))^{\times} \to ({\mathbb Z}/(p))^{\times}.\]
It is isomorphic to an additive group ${\mathbb Z}/(p^{r})$ and taking the inverse limit we have
\[\Gamma_{\infty}:=\lim_{\leftarrow}\Gamma_{r}\stackrel{\log}\simeq {\mathbb Z}_{p}.\]
Explicitly $\Gamma_{\infty}$ is the set of $p$-adic integers congruent to $1$ modulo $p$ and
\[\log x=\sum_{n=1}^{\infty}\frac{(-1)^{n-1}}{n}(x-1)^{n},\quad x\in \Gamma_{\infty}.\]
The isomorphism between $\Gamma_{r}$ and ${\mathbb Z}/(p^{r})$ is still denoted by $\log$.  There is an isomorphism
\begin{equation}{\mathbb Z}_{p}[\Gamma_{r}] \simeq {\mathbb Z}_{p}[t]/((1+t)^{p^r}-1),\end{equation}
 defined by
\[\varphi=\sum_{x\in\Gamma_{r}}\varphi(x)x\mapsto \sum_{x\in\Gamma_{r}}\varphi(x)(1+t)^{\log x}.\]
Let $\gamma\in \Gamma_{\infty}$ be a topological generator so that $\log\gamma=1$. Taking the inverse limit of (1) we have
\[{\mathbb Z}_{p}[[\Gamma_{\infty}]]:=\lim_{\leftarrow}{\mathbb Z}_{p}[\Gamma_{r}] \stackrel{\varpi_{p}}\simeq {\mathbb Z}_{p}[[t]], \quad \varpi_{p}(\gamma)=1+t.\]
Via $\varpi_{p}$ we sometimes identfy $t$ with $\gamma-1$. Putting 
\[t=e^{-s}-1=\sum_{n=1}^{\infty}\frac{(-s)^{n}}{n!},\]
 $\varpi_{p}$ yields an injective homomorphism of ${\mathbb C}_{p}$-algebras:
\begin{equation}\Lambda_{{\mathbb C}_{p}}:={\mathbb Z}_{p}[[\Gamma_{\infty}]]\otimes _{{\mathbb Z}_{p}}{\mathbb C}_{p}\stackrel{\tau_{p}}\hookrightarrow {\mathbb C}_{p}[[s]], \quad \tau_{p}(\gamma)=\sum_{n=0}^{\infty}\frac{(-s)^{n}}{n!}.\end{equation}
More explicitly using $p$-adic integral 
\begin{equation}\tau_{p}(\mu)(s)=\sum_{n=0}^{\infty}\frac{(-s)^{n}}{n!}\int_{\Gamma_{\infty}}(\log x)^{n}d\mu(x), \quad \mu\in \Lambda_{{\mathbb C}_{p}}.\end{equation}
 Let $a$ ($\neq 1$) be a $p$-adic integer that is congruent $1$ modulo $p$. We define {\it a Dirac measure} $\delta_{a}$ supported at $a$ to be
 \[\delta_{a}:=\lim_{\leftarrow}(\delta_{a})_{r}\in {\mathbb Z}_{p}[[\Gamma_{\infty}]], \quad (\delta_{a})_{r}=\sum_{x\in\Gamma_{r}}(\delta_{a})_{r}(x)x \in {\mathbb Z}_{p}[\Gamma_{r}],\]
 where $(\delta_{a})_{r}([a])=1$ and $(\delta_{a})_{r}(x)=0$ if $x\neq [a]$. Here $[a]$ is the image of $a$ by the natural 
 projection $\Gamma_{\infty}\to \Gamma_{r}$. A simple computation shows the following lemma.
\begin{lm}
\[\tau_{p}(\delta_{a})(s)=a^{-s}.\]
where
\[a^{-s}:=\sum_{n=0}^{\infty}\frac{(-s\log a)^{n}}{n!}\in {\mathbb C}_{p}[[s]].\]
\end{lm}
For an example let us take $a=\gamma$. Then (1) is a consequentce of (2).
Let ${\mathbb C}_{p}[\delta_{a}]$ be a subalgebra of  $\Lambda_{{\mathbb C}_{p}}$ generated by $\delta_{a}$, which is easily seen to be isomorphic to a ring of polynomials of one variable  whose coefficients are in ${\mathbb C}_{p}$. In fact let us consider a sub-semigroup $a^{{\mathbb N}}:=\{a^{m}\,|\,m\in {\mathbb Z},\,m\geq 0\}$ of $\Gamma_{\infty}$ which is  isomorphic to ${\mathbb N}:=\{x \in {\mathbb Z}|\,x \geq 0\}$. Then 
\[{\mathbb C}_{p}[\delta_{a}]= {\mathbb C}_{p}[a^{{\mathbb N}}]\subset \Lambda_{{\mathbb C}_{p}}.\]
We put $X_{a}=\delta_{a}-1$ and define
\[\frak{M}_p:=\tau_{p}^{-1}((s)), \quad \frak{N}:=\frak{M}_p \cap {\mathbb C}_{p}[X_{a}].\]
Then $\frak{M}_p$ and $\frak{N}$ are generated by $t$ and $X_{a}$, respectively. Let ${\mathbb C}_{p}[[X_{a}]]$ (resp. $\hat{\Lambda}_{{\mathbb C}_{p}}$) be a $\frak{N}$-adic (resp. $\frak{M}_p$-adic) completion of ${\mathbb C}_{p}[X_{a}]$ (resp. $\Lambda_{{\mathbb C}_{p}}$).  Note that 
\begin{equation}{\mathbb C}_{p}[X_{a}]/(X_{a}^{r})=\Lambda_{{\mathbb C}_{p}}/\frak{M}_p^{r}\stackrel{\tau_{p}}\simeq {\mathbb C}_{p}[s]/(s^{r}),\end{equation}
for every positive integer $r$. Passing to the inverse limit we see that 
\[{\mathbb C}_{p}[[X_{a}]]=\hat{\Lambda}_{{\mathbb C}_{p}},\]
and that $\tau_{p}$ is completed to an isomorphism:
\[{\mathbb C}_{p}[[X_{a}]]=\hat{\Lambda}_{{\mathbb C}_{p}}\stackrel{\hat{\tau}_{p}}\simeq  {\mathbb C}_{p}[[s]].\]
Summarizing we have proved the following.
\begin{prop} Let $a$ be an integer greater than $1$ that is congruent to $1$ modulo $p$. 
\begin{enumerate}
\item 
\[\hat{\Lambda}_{{\mathbb C}_{p}}={\mathbb C}_{p}[[X_{a}]], \quad X_{a}=\delta_{a}-1.\]
\item There is an injective homomorphism $\Lambda_{{\mathbb C}_{p}}\stackrel{\tau_{p}}\hookrightarrow {\mathbb C}_{p}[[s]]$ and it is completed to an isomorphism
\[{\mathbb C}_{p}[[X_{a}]]=\hat{\Lambda}_{{\mathbb C}_{p}}\stackrel{\hat{\tau}_{p}}\simeq  {\mathbb C}_{p}[[s]], \quad \hat{\tau}_{p}(X_{a})=a^{-s}-1.\]
In particular the natural map
\[\nu : \Lambda_{{\mathbb C}_{p}}\to \hat{\Lambda}_{{\mathbb C}_{p}},\]
is injective.
\end{enumerate}
\end{prop}
 Let ${\mathbb C}_p[[\Gamma_{\infty}]]$ be the inverse limit of a projective system $\{{\mathbb C}_p[\Gamma_{r}]\}_{r}$, which contains $\Lambda_{{\mathbb C}_{p}}$ as a subalgebra. Let $\alpha_{\Gamma_{r}}:{\mathbb C}_p[\Gamma_{r}]\to {\mathbb C}_p$ be the argumentation and 
\[\alpha_{\Gamma_{\infty}} : {\mathbb C}_p[[\Gamma_{\infty}]] \to {\mathbb C}_p,\]
be its inverse limit. Let $e_{0}: {\mathbb C}_p[[s]] \to {\mathbb C}_p$ be the evaluation at the origin: $e_{0}(f)=f(0)$. The following is derived from (3).
\begin{lm}
\[e_{0}(\tau_{p}(\mu))=\alpha_{\Gamma_{\infty}}(\mu),\quad \mu\in \Lambda_{{\mathbb C}_{p}}.\] 
\end{lm}
This shows that $\frak{M}_p$ is equal to the intersection of $\Lambda_{{\mathbb C}_p}$ with the argumentation ideal of ${\mathbb C}_{p}[[\Gamma_{\infty}]]$. 
We take a system of primitive $p^{n}$-th roots of the unit $\{\zeta_{p^n}\}_{n}$ satisfying $\zeta_{p^{n+1}}^{p}=\zeta_{p^n}$. According to a decomposition of Galois group:
\[{\rm Gal}({\mathbb Q}(\zeta_{p^{n+1}})/{\mathbb Q})\simeq ({\mathbb Z}/(p))^{\times}\times \Gamma_{n},\]
let ${\mathbb Q}_{n}$ be the abelian extension of ${\mathbb Q}$ contained in ${\mathbb Q}(\zeta_{p^{n+1}})$ such that ${\rm Gal}({\mathbb Q}_{n}/{\mathbb Q})\simeq \Gamma_{n}$ and ${\mathbb Q}_{\infty}$ their union: ${\mathbb Q}_{\infty}:=\cup_{n}{\mathbb Q}_{n}$. Then ${\rm Gal}({\mathbb Q}_{\infty}/{\mathbb Q})$ is isomorphic to $\Gamma_{\infty}$ and we will identify them. Then As we have explained before,
the projective limit ${\mathbb Z}_{p}[[{\rm Gal}({\mathbb Q}_{\infty}/{\mathbb Q})]]$ of $\{{\mathbb Z}_{p}[{\rm Gal}({\mathbb Q}_{n}/{\mathbb Q})]\}_{n}$
is isomorphic to ${\mathbb Z}_{p}[[t]]$. \\

Let $E$ be an elliptic curve defined over ${\mathbb Q}$ which has either good reduction or multiplicative reduction at $p$. Take a prime $l$ different from $p$ and let $\alpha_{E}\in {\mathbb Z}_{p}^{\times}$ and $\beta_{E}=p/{\alpha}_{E}\in p{\mathbb Z}_{p}$ be the eigenvalues of the $l$-adic representation of $p$-th power Frobenius on the Tate module $T_{l}(E)$ if $E$ has good ordinary reduction and $(\alpha_{E},\beta_{E})=(1,p)$ (resp. $(-1,-p)$) if $E$ has split (resp. non-split) multiplicative reduction at $p$. For a finite character $\chi$ of ${\rm Gal}({\mathbb Q}_{\infty}/{\mathbb Q})$ whose conductor $p^{n}$ let $W(\chi)$ be the Gauss sum:
\[W(\chi)=\sum_{\gamma\in {\rm Gal}({\mathbb Q}(\zeta_{p^{n}})/{\mathbb Q})}\chi(\gamma)\zeta_{p^{n}}^{\gamma}.\]
We fix an isomorphism $\sigma: {\mathbb C}_p \simeq {\mathbb C}$ such that
\[\sigma(z)=z,\quad z\in \bar{\mathbb Q}.\]
\begin{fact}(\cite{MTT}) There is the unique element $\mu_{E,p}$ of ${\mathbb Z}_{p}[[{\rm Gal}({\mathbb Q}_{\infty}/{\mathbb Q})]]\otimes_{{\mathbb Z}_{p}}{\mathbb Q}_{p}$ satisfying the following properties:
\begin{enumerate}
\item If $E$ has good ordinary reduction at $p$,
\[\sigma({\mathbf 1}(\mu_{E,p}))=(1-\alpha_{E}^{-1})^{2}\frac{L(E,1)}{\Omega_{E}}.\]

\item If $E$ has multiplicative reduction at $p$,
\[\sigma({\mathbf 1}(\mu_{E,p}))=(1-\alpha_{E}^{-1})\frac{L(E,1)}{\Omega_{E}}.\]

\item Let $\chi$ be a character  of ${\rm Gal}({\mathbb Q}_{\infty}/{\mathbb Q})$ of finite order whose conductor $p^{n}>1$. Then
\[\sigma(\chi(\mu_{E,p}))=\frac{W(\chi)}{\alpha_{E}^{n}}\frac{L(E,\chi^{-1},1)}{\Omega_{E}}.\]
\end{enumerate}
\end{fact}
Here ${\mathbf 1}$ is the trivial character and $\Omega_{E}$ is the fundamental real period of $E$. The $p$-adic L-function ${\mathcal L}_{E,p}$ of $E$ is defined to be
\[{\mathcal L}_{E,p}:=\tau_{p}(\mu_{E,p})\in {\mathbb C}_{p}[[s]].\]
Let $\phi$ and $\iota$ be automorphisms of $\Gamma_{\infty}$ defined to be
\[\phi(x)=x^{p-1},\quad \iota(x)=x^{-1},\quad x\in \Gamma_{\infty},\]
 and the induced automorphisms on $\Lambda_{{\mathbb C}_p}$ are still denoted by the same character.  
\begin{lm} Let $\mu\in\Lambda_{{\mathbb C}_p}$.
\begin{enumerate}
\item \[\chi(\phi(\mu))=\chi^{p-1}(\mu), \quad \chi(\iota(\mu))=\chi^{-1}(\mu).\]
\item \[\tau_{p}(\phi(\mu))(s)=\tau_{p}(\mu)((p-1)s),\quad \tau_{p}(\iota(\mu))(s)=\tau_{p}(\mu)(-s).\]
\end{enumerate}
\end{lm}
The following follows from {\bf Fact 2.1} and {\bf Lemma 2.3}
\begin{prop}
\begin{enumerate}
\item If $E$ has good ordinary reduction at $p$,
\[\sigma({\mathbf 1}(\iota\phi(\mu_{E,p})))=(1-\alpha_{E}^{-1})^{2}\frac{L(E,1)}{\Omega_{E}}.\]

\item If $E$ has multiplicative reduction at $p$,
\[\sigma({\mathbf 1}(\iota\phi(\mu_{E,p})))=(1-\alpha_{E}^{-1})\frac{L(E,1)}{\Omega_{E}}.\]
\item Let $\chi$ be a character of ${\rm Gal}({\mathbb Q}_{\infty}/{\mathbb Q})$ of finite order whose conductor $p^{n}>1$. Then
\[\sigma(\chi(\iota\phi(\mu_{E,p})))=\frac{W(\chi^{1-p})}{\alpha_{E}^{n}}\frac{L(E,\chi^{p-1},1)}{\Omega_{E}}.\]
\item
\[\tau_{p}(\iota\phi(\mu_{E,p}))(s)={\mathcal L}_{E,p}((1-p)s).\]
\end{enumerate}
\end{prop}
\subsection{An archimedian analog of the $p$-adic measure}
We want to develop an analog of the $p$-adic measure theory over ${\mathbb C}$ for a certain Dirichlet series. For $\rho\in{\mathbb R}$ let $H_{\rho}$ be a right half plane defined by
\[H_{\rho}:=\{z\in{\mathbb C}\,:\, {\rm Re}\,z > \rho \}.\]
\begin{df} A Dirichlet series
\[A(z)=\sum_{n=1}^{\infty}a_{n}n^{-z},\quad a_{n}\in{\mathbb C},\]
will be called {\rm regular} if it satisfies the following conditions:
\begin{enumerate}
\item  $a_{n}=0$ if $p$ divides $n$.
\item There is a positive real number $\rho$ (which may depend on $A(z)$) such that $A(z)$ absolutely converges on $H_{\rho}$ and is analytically continued to the whole plane as an entire function.
\item For every positive integer $r$ and $a$ such that $1\leq a \leq p^{r}$ 
\[A_{r,a}(z):=\sum_{k=1, k\equiv a(p^{r})}^{\infty}a_{k}k^{-z}\]
is also continued to the whole plane as an entire function.
\end{enumerate}
\end{df}
We will denote the set of regular Dirichlet series by ${\mathcal R}$. By definition it is contained in the commutative algebra ${\mathcal O}_{\mathbb C}$ of holomorphic functions on ${\mathbb C}$. In fact it is a subalgebra. We only check that it is closed by a multiplication. Let $A(z)=\sum_{l=1}^{\infty}a_{l}l^{-z}$ and $B(z)=\sum_{m=1}^{\infty}b_{m}m^{-z}$ be regular Dirichlet series and $C(z)=\sum_{n=1}^{\infty}c_{n}n^{-z}$ their product. It is obvious that $C(z)$ satisfies (1) and (2). For $1\leq c \leq p^{r}$ a simple computation shows
\begin{equation}C_{r,c}(z)=\sum_{ab\equiv c (p^r), 1\leq a,b \leq p^r}A_{r,a}(z)B_{r,b}(z),\end{equation}
and this implies (3). For $A(z)=\sum_{n=1}^{\infty}a_{n}n^{-z} \in{\mathcal R}$, we define
\[\mu_{{\mathcal O},r}(A):=\sum_{a=1}^{p^{r}}A_{r,a}(z)[a]_{r} \in {\mathcal O}_{\mathbb C}[({\mathbb Z}/(p^{r}))^{\times}],\]
where $[\cdot]_{r}$ represents the residue class.
Then $\{\mu_{{\mathcal O},r}(A)\}_{r}$ forms a projective system and we set
\[\mu_{{\mathcal O}}(A):=\lim_{\leftarrow}\mu_{{\mathcal O},r}(A) \in {\mathcal O}_{\mathbb C}[[{\mathbb Z}_{p}^{\times}]].\]
Thus we have a map
\[\mu_{{\mathcal O}} : {\mathcal R} \to {\mathcal O}_{\mathbb C}[[{\mathbb Z}_{p}^{\times}]],\]
which is a homomorphism of algebras by (5). Let $\tilde{\Lambda}_{{\mathcal O}}$ be its image. Note that a regular Dirichlet series $A$ is recovered from $\mu_{{\mathcal O}}(A)$. In fact we associate a function
$\sum_{a=1}^{p^{r}}A_{r,a}(z)a^{-s}$ to $\mu_{{\mathcal O},r}(A)$, that is $\sum_{a=1}^{p^{r}}\sum_{k=1,k\equiv a(p^{r})}a_{k}k^{-z}a^{-s}$ on $H_{\rho}\times H_{\rho}$ ($\rho$ is sufficiently large). Therefore 
\[\lim_{r\to \infty}\sum_{a=1}^{p^{r}}\sum_{k=1,k\equiv a(p^{r})}a_{k}k^{-z}a^{-s}=\sum_{n=1,(n,p)=1}^{\infty}a_{n}n^{-(s+z)}=A(s+z),\]
on $H_{\rho}\times H_{\rho}$ and by analytic continuation we have a map
\begin{equation}\tilde{\tau}_{\mathcal O} : \tilde{\Lambda}_{{\mathcal O}} \to {\mathcal O}_{{\mathbb C}\times{\mathbb C}},\quad \tilde{\tau}_{\mathcal O}(\mu_{{\mathcal O}}(A))=A(s+z),\end{equation}
where ${\mathcal O}_{{\mathbb C}\times{\mathbb C}}$ is the set of holomorphic functions on ${\mathbb C}\times{\mathbb C}$. Thus 
\[\tilde{\tau}_{\mathcal O}\mu_{\mathcal O}=p^{*},\]
where
\[p : {\mathbb C}\times{\mathbb C} \to {\mathbb C}, \quad p(s,z)=s+z,\]
and since $p^{*}$ is an injective homomorphism, ${\mathcal R}\stackrel{\mu_{\mathcal O}}\to \tilde{\Lambda}_{{\mathcal O}}$ is an isomorphism and $\tilde{\tau}_{\mathcal O}$ is injective. Let
\[e : {\mathcal O}_{\mathbb C}[[{\mathbb Z}_{p}^{\times}]] \to {\mathbb C}[[{\mathbb Z}_{p}^{\times}]],\]
be a homomorphism induced by
\[{\mathcal O}_{\mathbb C} \to {\mathbb C}, \quad f\mapsto f(1),\]
and set $\tilde{\tau}_{\mathbb C}:=\tilde{\tau}_{\mathcal O}\otimes_{{\mathcal O},e}{\mathbb C}$. Let $\mu_{\mathbb C}$ be the composition of \[{\mathcal R}\stackrel{\mu_{{\mathcal O}}}\simeq \tilde{\Lambda}_{\mathcal O} \stackrel{e}\to \tilde{\Lambda}_{\mathbb C}.\]
Then by (6)
\[\tilde{\tau}_{\mathbb C}\circ \mu_{\mathbb C}(A)(s)=A(1+s),\quad A\in {\mathcal R},\]
and we have proved the following proposition.

\begin{prop}
\[\mu_{\mathbb C} : {\mathcal R} \to \tilde{\Lambda}_{\mathbb C},\]
and
\[\tilde{\tau}_{\mathbb C} : \tilde{\Lambda}_{\mathbb C} \to {\mathcal O}_{\mathbb C},\]
are an isomorphism and an injective homomorphism, respectively. Moreover
\[\tilde{\tau}_{\mathbb C}\mu_{\mathbb C}(A)=A(1+s),\quad A\in{\mathcal R}.\]
\end{prop}
Here is an example. Let $b (\neq 1)$ be a positive integer prime to $p$ and set $D_{b}:=b^{1-z}$,  which is a regular Dirichlet series. Then
\begin{equation}\mu_{\mathbb C}(D_{b})=\delta_{b},\quad \tilde{\tau}_{\mathbb C}(D_{b})=b^{-s}.\end{equation}
Let $\tilde{\tau}_{\infty} : \tilde{\Lambda}_{\mathbb C} \to {\mathbb C}[[s]]$ be the composition of $\tilde{\tau}_{\mathbb C}$ with the Taylor expansion at the origin:
\[{\mathcal O}_{{\mathbb C}} \to {\mathbb C}[[s]],\quad f \mapsto \sum_{n=0}^{\infty}a_{n}(f)s^{n},\]
which is injective by {\bf Proposition 2.3}. It is easy to check that
\[\tilde{\tau}_{\infty}(A)(0)=\alpha_{\tilde{\Lambda}_{\mathbb C}}(A),\quad A\in \tilde{\Lambda}_{\mathbb C},\]
where $\alpha_{\tilde{\Lambda}_{\mathbb C}}$ is the restriction of the argumentation of ${\mathbb C}[[{\mathbb Z}_{p}^{\times}]]$. In particular we see that the kernel $\frak{M}_{\tilde{\Lambda}_{\mathbb C}}$ of $\alpha_{\tilde{\Lambda}_{\mathbb C}}$ is equal to $\tilde{\tau}_{\infty}^{-1}(s)$. Take a positive integer $b$ as above
and set $Y_{b}:=\delta_{b}-1$. Then as we have seen in the previous section ${\mathbb C}[Y_{b}]$ is a subalgebra of $\tilde{\Lambda}_{\mathbb C}$ which is isomorphic to the polynomial ring of one variable and it is esy to see that
\[(Y_{b})={\mathbb C}[Y_{b}]\cap \frak{M}_{\tilde{\Lambda}_{\mathbb C}}.\]
 Thus we have shown that the each arrow of
\begin{equation}{\mathbb C}[Y_{b}]/ (Y_{b}^{r}) \to \tilde{\Lambda}_{\mathbb C}/\frak{M}_{\tilde{\Lambda}_{\mathbb C}}^{r}  \stackrel{\tilde{\tau}_{\infty}} \to {\mathbb C}[[s]]/(s^{r})\end{equation}
is injective. Since the dimension of both side are equal they are isomorphic. The following is an archimedian analog of {\bf Proposition 2.1}.
\begin{prop}
\begin{enumerate}
\item \[ {\mathbb C}[[Y_{b}]]=\hat{\tilde{\Lambda}}_{\mathbb C}\stackrel{\hat{\tilde{\tau}}_{\infty}}\simeq {\mathbb C}[[s]].\]
\item The natural map $\tilde{\Lambda}_{\mathbb C} \to \hat{\tilde{\Lambda}}_{\mathbb C}$ is injective.
\end{enumerate}
\end{prop} 
{\bf Proof.} Take the limit of (8) we obtain (1).  
\[\begin{CD}
 {\mathcal R}@>\text{$\tilde{\tau}_{\infty}\circ \mu_{\mathbb C}$}>> {\mathbb C}[[s]]\\
 @V\text{$\mu_{\mathbb C}$}VV  @V\text{$\hat{\tilde{\tau}}_{\infty}^{-1}$}VV\\
   \hat{{\Lambda}}_{\mathbb C}@>>>  \hat{\tilde{\Lambda}}_{\mathbb C},
 \end{CD}\]
Notice that $\tilde{\tau}_{\infty}\circ \mu_{\mathbb C} : {\mathcal R} \to {\mathbb C}[[s]]$ is the Taylor expansion at $s=1$ and it is injective. Since vertical arrows are isomorphisms we obtain (2).
\begin{flushright}
$\Box$
\end{flushright}
But there is a slight difference. Namely ${\Lambda}_{{\mathbb C}_{p}}$ is a subalgebra of ${\mathbb C}_{p}[[\Gamma_{\infty}]]$ but ${\tilde{\Lambda}}_{\mathbb C}$ is contained in ${\mathbb C}[[{\mathbb Z}_{p}^{\times}]]$.
Let
\[\phi : {\mathbb C}[[{\mathbb Z}_{p}^{\times}]] \to {\mathbb C}[[\Gamma_{\infty}]],\]
be a homomorphism induced by 
\[\phi : {\mathbb Z}_{p}^{\times} \to \Gamma_{\infty}\quad \phi(x)=x^{p-1},\]
and we define
\[\Lambda_{\mathbb C}:=\phi(\tilde{\Lambda}_{\mathbb C}).\]
Let $\alpha_{\Lambda_{\mathbb C}}$ be the restriction of the argumentation of ${\mathbb C}[[\Gamma_{\infty}]]$ to $\Lambda_{\mathbb C}$ and $\frak{M}_{\infty}$ its kernel. Set $a:=b^{p-1}$ and $X_{a}:=\delta_{a}-1$. Then ${\mathbb C}[X_{a}]$ is a subalgebra of $\Lambda_{\mathbb C}$ isomorphic to a polynomial ring of one variable and $\phi$ yields a isomorphism:
\[\phi : {\mathbb C}[Y_{b}] \to  {\mathbb C}[X_{a}], \quad \phi(Y_{b})=X_{a}.\]
Since 
\[(X_{a})={\mathbb C}[X_{a}]\cap \frak{M}_{\infty}, \quad \phi^{-1}(\frak{M}_{\infty})=\frak{M}_{\tilde{\Lambda}_{\mathbb C}},\]
we obtain the following diagram.
\begin{equation}\begin{CD}
 {\mathbb C}[Y_{b}]/ (Y_{b}^{r})@>\text{$\tilde{i}$}>> \tilde{\Lambda}_{\mathbb C}/\frak{M}_{\tilde{\Lambda}_{\mathbb C}}^{r} @> \text{$\tilde{\tau}_{\infty}$}>> {\mathbb C}[[s]]/(s^{r})\\
 @V\text{$\phi$}VV  @V\text{$\phi$}VV\\
  {\mathbb C}[X_{a}]/ (X_{a}^{r}) @>\text{$i$}>>  {\Lambda}_{\mathbb C}/\frak{M}_{\infty}^{r}
 \end{CD}\end{equation}
 By (8) upper horizontal arrows are isomorphisms. Since the both $\phi$ are isomorphic so is  $i$. (The reason that the right $\phi$ is isomorphic is $\phi : \tilde{\Lambda}_{\mathbb C} \to {\Lambda}_{\mathbb C}$ is surjective and $\phi^{-1}(\frak{M}_{\infty})=\frak{M}_{\tilde{\Lambda}_{\mathbb C}}$.) Define $\hat{\tau}_{\infty}$ and $\phi$ to be 
\[\hat{\tau}_{\infty} : {\mathbb C}[[X_{a}]] \to {\mathbb C}[[s]], \quad \hat{\tau}_{\infty}(X_{a})=a^{-s}-1,\]
and
\[\phi : {\mathbb C}[[s]] \to {\mathbb C}[[s]], \quad \phi(f)(s)=f((p-1)s).\]
Take the limit of (9) and we have
\[\begin{CD}
{\tilde{\Lambda}}_{\mathbb C}@>\text{$\tilde{\nu}$}>> {\mathbb C}[[Y_{b}]]=\hat{\tilde{\Lambda}}_{\mathbb C}@>\text{$\hat{\tilde{\tau}}_{\infty}$}>> {\mathbb C}[[s]]\\
 @V\text{$\phi$}VV @V\text{$\phi$}VV  @V\text{$\phi$}VV\\
{{\Lambda}}_{\mathbb C} @>\text{$\nu$}>>  {\mathbb C}[[X_{a}]]=\hat{{\Lambda}}_{\mathbb C}@>\text{$\hat{{\tau}}_{\infty}$}>>  {\mathbb C}[[s]],
 \end{CD}\]
 where $\tilde{\nu}$ and $\nu$ are natural homomorphisms and every arrow in the right rectangle is an isomorphism. As we have shown in {\bf Proposition 2.4}  $\tilde{\nu}$ is injective and by definition the most left $\phi$ is surjective. A diagram chasing shows that $\nu$ is injective and that the most left $\phi$ is also an isomorphism.  Let us denote ${\tau}_{\infty}=\hat{{\tau}}_{\infty}\circ \nu$. Here is an archimedian analog of {\bf Proposition 2.2}.
 \begin{prop}
 \begin{enumerate}
 \item
 \[{\mathbb C}[[X_{a}]]=\hat{{\Lambda}}_{\mathbb C} \stackrel{\hat{{\tau}}_{\infty}} \simeq {\mathbb C}[[s]], \quad  \hat{\tau}_{\infty}(X_{a})=a^{-s}-1.\]
\item The natural map $\nu : {{\Lambda}}_{\mathbb C} \to \hat{{\Lambda}}_{\mathbb C}$ is injective.
\item For $A\in{\mathcal R}$ ${\tau}_{\infty}(\phi\circ \mu_{\mathbb C}(A))$ is the Taylor expansion of $A(1+(p-1)s)$ at the origin and
\[\alpha_{\Lambda_{\mathbb C}}(\phi\circ \mu_{\mathbb C}(A))=A(1).\]

 \end{enumerate}
 \end{prop}
 \subsection{An archimedian measure of a cusp form}
 Let $f=\sum_{n=1}^{\infty}a_n(f)q^n$ be a cusp form of weight 2 and level $N$. Let us fix a positive integer $m$ and for an integer $0\leq a\leq m-1$ we put
 \[\phi_m^a=f(z+\frac{a}{m}), \quad f_m^a=\sum_{n=1, n\equiv a(m)}^{\infty}a_n(f)q^n.\]
 \begin{lm}
 $\{\phi_m^0,\cdots, \phi_m^{m-1}\}$ and $\{f_m^0,\cdots, f_m^{m-1}\}$ span the same vector space over ${\mathbb Q}(\zeta_m)$.
 \end{lm}
 {\bf Proof.} Putting $\zeta_m=\exp(2\pi i/m)$ a simple computation shows
 \[\phi_m^a=\sum_{n=1}^{\infty}a_n(f)q^{n}\zeta_{m}^{an}=\sum_{k=0}^{m-1}\zeta_m^{ak}f_m^k.\]
 Therefore
 \[ \left(\begin{array}{c}\phi_m^{0} \\ \vdots \\ \phi_{m}^{m-1}\end{array}\right)=A\left(\begin{array}{c}f_m^{0} \\ \vdots \\ f_{m}^{m-1}\end{array}\right)\]
 where 
 \[A=\left(\begin{array}{cccc}1 & 1 & \cdots & 1 \\1 & \zeta_m & \cdots & \zeta_m^{m-1} \\ \vdots & \vdots & \ddots & \vdots \\1 & \zeta_m^{m-1} & \cdots & \zeta_m^{(m-1)^2}\end{array}\right).\]
 Since $A$ is regular, we obtain the claim. $\Box$\\
 
Let $E$ be an elliptic curve defined over ${\mathbb Q}$ and $N_{E}$ its conductor. Since it is modular \cite{BCDT} there is a cusp form $f_{E}$ of weight $2$ and level $N_{E}$ associated to $E$ and let 
\[f_{E}=\sum_{n=1}^{\infty}a_{n}(E)q^{n}, \quad q=\exp(2\pi iz)\]
 be the Fourier expansion at $i\infty$.  Removing the Euler factor $L_{p}(E,z)$ at $p$ we define {\it the modified $L$-function} of $E$ to be
\[L^{\dagger}(E,z)=\sum_{n=1, (n,p)=1}^{\infty}a_{n}(E)n^{-z}.\]
Since its partial Dirichlet series
\[(L^{\dagger}(E,z))_{r,a}:=\sum_{k=1,k\equiv a (p^{r})}^{\infty}a_{k}(E)k^{-z}\]
is a Mellin transform of $(f_{E})_{p^{r}}^{a}$ {\bf Lemma 2.4} and the cuspidality of $f_{E}$ imply that $L^{\dagger}(E,z)$ is a regular Dirichlet series.  Now we set
\[\mu_{E,\infty}:=\mu_{\mathbb C}(L^{\dagger}(E,z))\in\tilde{\Lambda}_{\mathbb C}.\]
Let $\chi$ be a character of ${\mathbb Z}_{p}^{\times}$ of finite order whose conductor is $p^{r}$. It defines a homomorphism
\[ \chi : {\mathbb C}[({\mathbb Z}/(p^{r}))^{\times}] \to {\mathbb C},\]
and the composition it with the projection ${\mathbb C}[[{\mathbb Z}_{p}^{\times}]] \to {\mathbb C}[({\mathbb Z}/(p^{r}))^{\times}]$ is still denoted by $\chi$. Then 
\[\chi(\mu_{E,\infty})=L(E,\chi,1),\]
for a non-trivial $\chi$ and  
\[{\bf 1}(\mu_{E,\infty})=\alpha_{\Lambda_{\mathbb C}}(\mu_{E,\infty})=L^{\dagger}(E,1).\]
Thus we see
\[\chi(\phi(\mu_{E,\infty}))=L(E,\chi^{p-1},1),\]
by {\bf Lemma 2.3}. {\bf Proposition 2.5} implies the following theorem, which should be compared to {\bf Fact 2.1}.

\begin{thm} Let $E$ be an elliptic curve defined over ${\mathbb Q}$. Then $\phi(\mu_{E,\infty})\in \Lambda_{\mathbb C}$ satisfies the following properties.
\begin{enumerate}
\item 
\[\chi(\phi(\mu_{E,\infty}))=L(E,\chi^{p-1},1),\]
for a non-trivial character $\chi$ of $\Gamma_{\infty}$ of finite order.
\item 
\[{\bf 1}(\phi(\mu_{E,\infty}))=L^{\dagger}(E,1)=\frac{L(E,1)}{L_{p}(E,1)}.\]
\item 
\[\tau_{\infty}(\phi(\mu_{E,\infty}))(s)=L^{\dagger}(E,1+(p-1)s).\]
\end{enumerate}
\end{thm}
Our measure is related to Kato's system which will be recalled below.
Let $T_p(E)$ and $V_p(E)$ be the Tate module of $E$:
\[T_{p}(E)=\lim_{\leftarrow}E[p^n], \quad V_{p}(E)=T_{p}(E)\otimes_{{\mathbb Z}_p}{\mathbb Q}_p.\]
${\mathbb Q}_{n,p}$ denotes the completion of ${\mathbb Q}_{n}$ by the unique prime on $p$. Let $H^{1}_{S}({\mathbb Q}_{n,p}, V_{p}(E))\subset H^{1}({\mathbb Q}_{n,p}, V_{p}(E))$ be the image of $E({\mathbb Q}_{n,p})\otimes {\mathbb Q}_p$ by the Kummer map and $H^{1}_{S}({\mathbb Q}_{n,p}, T_{p}(E))$ its intersection of with $H^{1}({\mathbb Q}_{n,p}, T_{p}(E))$. For $M=T_{p}(E)$ or $V_{p}(E)$ we define
\[H^{1}_{\slash{S}}({\mathbb Q}_{n,p}, M)=H^{1}({\mathbb Q}_{n,p}, M)\slash H^{1}_{S}({\mathbb Q}_{n,p}, M).\]
Let $\omega_E$ be the canonical invariant differential associated to the minimal Weierstrauss model of $E$. Then there is an isomorphism called {\it dual exponential map}:
\[ \exp^{*} : H^{1}_{\slash{S}}({\mathbb Q}_{n,p}, V_{p}(E)) \tilde{\to} {\mathbb Q}_{n,p}\omega_E.\]
Restrict the composition of $\exp^{*}$ with 
\[{\mathbb Q}_{n,p}\omega_E \tilde{\to} {\mathbb Q}_{n,p},\quad a \omega_{E} \mapsto a\]
to $H^{1}_{\slash{S}}({\mathbb Q}_{n,p}, T_{p}(E))$
and we obtain a map
\[\exp^{*}_{\omega_E} : H^{1}_{\slash{S}}({\mathbb Q}_{n,p}, T_{p}(E)) \to {\mathbb Q}_{n,p}.\]
\begin{fact}(\cite{Kato}, \cite{Rubin1998} {\bf Corollary 7.2})
For every $n$ there is $c_n\in H^1({\mathbb Q}_n, T_{p}(E))$ satisfying following properties.
\begin{enumerate}
\item 
\[{\rm Cor}_{n,n+1}(c_{n+1})=c_n.\]
where ${\rm Cor}_{n,n+1} : H^1({\mathbb Q}_{n+1}, T_{p}(E)) \to H^1({\mathbb Q}_n, T_{p}(E))$ is a corestriction map.
\item Let ${\rm loc}_{p}^{ram}$ be the composition of
\[H^1({\mathbb Q}_n, T_{p}(E)) \stackrel{{\rm loc}_{p}}\to H^1({\mathbb Q}_{n,p}, T_{p}(E)) \to H^{1}_{\slash{S}}({\mathbb Q}_{n,p}, T_{p}(E)),\]
where the first arrow is the localization and the second is the natural projection. Then
\[\sigma(\sum_{\gamma\in \Gamma_n}\chi(\gamma)\exp^{*}_{\omega_E} ({\rm loc}_{p}^{ram}(c_n^{\gamma})))=\frac{r_E}{\Omega_E}L_{(pN_E)}(E,\chi,1),\]
 for any character $\chi$ of $\Gamma_n$. Here $r_E$ is a positive integer which depends only on $E$ and $\Omega_E$ is the fundamental real period. $ L_{(pN_E)}(E,\chi,s)$ is a function obtained from $L(E,\chi,s)$  removing Euler factors at primes which divide $pN_E$.
\end{enumerate}
\end{fact}
Let 
\[\kappa_n=\sum_{\gamma\in \Gamma_n}\exp^{*}_{\omega_E} ({\rm loc}_{p}^{ram}(c_n^{\gamma}))\gamma \in {\mathbb C}_{p}[\Gamma_n].\]
Then $\{\kappa_n\}_n$ forms a projective system by {\bf Fact 2.2} and \[\chi(\sigma(\phi(\kappa_{\infty})))=\frac{r_E}{\Omega_E}L_{(pN_E)}(E,\chi^{p-1},1),\quad \kappa_{\infty}=\lim_{\leftarrow}\kappa_{n},\]
for any finite character $\chi$ of $\Gamma_{\infty}$ by {\bf Lemma 2.3}.
Since
\[L_{(pN_E)}(E,\chi,1)=\prod_{q|pN_E}P_q(q^{-1}\chi({\rm Fr}_q))L(E,\chi,1),\]
where
\[P_q(t)=
 \begin{cases}
 1-t, & \text{if $E$ has a split multiplicative reduction at $q$:}\\
 1+t, & \text{if $E$ has a non-split multiplicative reduction at $q$:}\\
 1, & \text{if $E$ has an additive reduction at $q$,}
 \end{cases}\]
{\bf Theorem 2.1} implies the following result.
\begin{prop}
\[\sigma(\phi(\kappa_{\infty}))=\frac{r_E}{\Omega_E}\prod_{q|N_E, q\neq p}\phi(P_q(q^{-1}{\rm Fr}_q)\cdot \mu_{E,\infty}).\]
\end{prop}

\subsection{A motivation of the conjecture}
Now we are ready to explain a motivation of {\bf Conjecture 1.1}. We will fix   a positive integer $b$ greater than 1 which is prime to $p$ and set $a=b^{p-1}$. In  order to make a distinction between $p$-adic and archimedian logarithm of $a$ we denote them by $l_{p}(a)\in {\mathbb C}_p$ and $l_{\infty}(a)\in {\mathbb C}$, respectively. We refer an isomorphism $\sigma : {\mathbb C}_p \to 
{\mathbb C}$ is {\it normalized} if it satisfies
\begin{enumerate}
\item \[\sigma(z)=z,\quad z\in \bar{\mathbb Q}.\]
\item \[\sigma(l_{p}(a))=l_{\infty}(a).\]
\end{enumerate}
Let us extends $\sigma$ to  ${\mathbb C}_p[[s]]$ and ${\mathbb C}_p[[X_{a}]]$ by
\[\sigma(\sum_{n}c_{n}s^{n})=\sum_{n}(c_n)s^{n},\]
and
\[\sigma(\sum_{n}d_{n}X_{a}^{n})=\sum_{n}\sigma(d_n)X_{a}^{n}.\]
Then (2) implies
\[\sigma(a^{-s})=a^{-s},\]
as a power series of $s$ and in particular
\begin{equation}\sigma\circ \hat{\tau}_{p}=\hat{\tau}_{\infty}\circ \sigma.\end{equation}
Finally we define
\[\sigma : {\mathbb C}_p[[\Gamma_{\infty}]] \to {\mathbb C}[[\Gamma_{\infty}]]\]
to be the inverse limit of 
\[ \sigma : {\mathbb C}_p[\Gamma_{r}] \to {\mathbb C}[\Gamma_{r}],\quad 
\sigma(\sum_{\gamma\in \Gamma_{r}}c_{\gamma}\gamma)=\sum _{\gamma\in \Gamma_{r}}\sigma(c_{\gamma})\gamma.\] 
Since $\sigma$ is normalized,
\[\sigma(\chi(c))=\chi(\sigma(c)),\quad c\in {\mathbb C}_p[[\Gamma_{\infty}]],\]
for a finite character $\chi$ of $\Gamma_{\infty}$. Let $\Lambda$ be a ${\mathbb C}$-subalgebra of ${\mathbb C}[[\Gamma_{\infty}]]$ generated by
$\sigma(\Lambda_{{\mathbb C}_p})$ and $\Lambda_{\mathbb C}$ and $\frak{M}_{\Lambda}$ the intersection of it and the argumentation ideal of ${\mathbb C}[[\Gamma_{\infty}]]$. Then $\Lambda$ contains ${\mathbb C}[X_{a}]$ and 
\begin{equation}\Lambda/\frak{M}_{\Lambda}^{r}={\mathbb C}[X_{a}]/(X_{a})^{r} \stackrel {\hat{\tau}_{\infty}} \simeq {\mathbb C}[[s]]/(s^{r}),\quad \forall r\end{equation}
by {\bf Proposition 2.1} and {\bf Proposition 2.4}. Therefore we see that, for $\lambda\in\Lambda$, 
\begin{equation}{\rm Min}\{k\,|\, \lambda\in \frak{M}_{\Lambda}^{k}\}={\rm ord}_{s=0}\hat{\tau}_{\infty}(\lambda).\end{equation}
Let $E$ and $E^{\prime}$ be elliptic curves satisfying the assumption of {\bf Conjecture 1.1}. By {\bf Proposition 2.2} and {\bf Theorem 2.1} we see that
\begin{equation}\frac{\phi(\mu_{E^{\prime},\infty})}{\Omega_{E^{\prime}}}\cdot \sigma(\iota\phi(\mu_{E,p}))=\frac{\phi(\mu_{E,\infty})}{\Omega_{E}}\cdot \sigma(\iota\phi(\mu_{E^{\prime},p}))\in \Lambda.\end{equation}
Now {\bf Conjecutre 1.1} will be derived by (10), (12), {\bf Proposition 2.2} and {\bf Theorem 2.1}.\\

 Unfortunately there seems to be a mistake in the above argument. In fact suppose that both $E$ and $E^{\prime}$ have a split multiplicative reduction at $p$ and that their L-function does not vanish at $s=1$. Take the image of (13) by $\hat{\tau}_{\infty}$ and the formula of the first derivative of $p$-adic L-function (\cite{Greenberg-Stevens}, \cite{Kobayashi2006}) will tell us that their L-invariants should be equal! Since we cannot solve this puzzle the above argument should be considered as only an explanation and not a proof. But suppose that $E^{\prime}$ is a quadratic twist of $E$ (as we will see in the next section, in order to derive consequences from {\bf Conjecture 1.1}, we may impose this). Because they are isomorphic over $\overline{\mathbb Q}$ their $j$-invariants are equal and  so are L-invariants since the Tate period is determined by the $j$-invariant (\cite{MTT} {\bf Chapter II} \$1 (2)).
 \begin{conjecture}
Let $E$ and $E^{\prime}$ be elliptic curves defined over ${\mathbb Q}$. Suppose that $E^{\prime}$ is a quadratic twist of $E$ and that they have ordinary reduction of the same type at $p$. Then
\[\Omega_{E^{\prime}}L(E,1-s)\sigma({\mathcal L}_{E^{\prime},p}(s))=\Omega_{E}L(E^{\prime}, 1-s)\sigma({\mathcal L}_{E,p}(s)).\]
\end{conjecture}
Although our argument may be incomplete it will explain why an extra zero appears. We will follow the notation of \cite{Rubin1998} {\bf Appendix}. Suppose that $E$ has a split multiplicative reduction at $p$. For simplicity  we will omit to write the isomorphism $\sigma : {\mathbb C}_p \to {\mathbb C}$ and will make no distinction between $\hat{\tau}_{p}$ and $\hat{\tau}_{\infty}$, which will be denoted by $\hat{\tau}$.
 Let us put
 \[x_{n}={\rm Tr}_{{\mathbb Q}(\zeta_{p^{n+1}})/{\mathbb Q}_{n,p}}(\sum_{k=0}^{n}\frac{\zeta_{p^{n+1-k}}-1}{p^k}+\frac{p}{p-1})\in {\mathbb Q}_{n,p}\]
 and 
 \[w_{n}=\sum_{\gamma\in \Gamma_{n}}x_{n}^{\gamma}\gamma \in {\mathbb Q}_{n,p}[\Gamma_{n}].\]
 One may check that $\{x_{n}\}_{n}$ is compatible with the corestrictions and that $\{w_{n}\}_{n}$ forms a projective system. According to Coleman define
 \[{\rm Col}_{n}(\kappa_{n}):=w_{n}\iota(\kappa_{n})=\sum_{\gamma\in\Gamma_{n}}({\rm Tr}_{{\mathbb Q}_{n,p}/{\mathbb Q}_{p}}x_{n}^{\gamma})\exp^{*}_{\omega_E} ({\rm loc}_{p}^{ram}(c_n))\gamma^{-1}.\]
 It is known that it is contained in ${\mathbb Z}_{p}[\Gamma_{n}]$ and
 that (\cite{Rubin1998} {\bf Corollary 7.2})
 \[r_{E}\prod_{q|N_E, q\neq p}P_q(q^{s-1}){\mathcal L}_{E,p}(s)=\hat{\tau}({\rm Col}_{\infty}(\kappa_{\infty})), \quad {\rm Col}_{\infty}(\kappa_{\infty})=\lim_{\leftarrow}{\rm Col}_{n}(\kappa_{n}).\]
 Therefore
 \[r_{E}\prod_{q|N_E, q\neq p}P_q(q^{(p-1)s-1}){\mathcal L}_{E,p}((p-1)s)=\hat{\tau}\phi({\rm Col}_{\infty}(\kappa_{\infty})).\]
 Since ${\rm Col}_{\infty}(\kappa_{\infty})=w_{\infty} \iota(\kappa_{\infty})$ we formally obtain
 \[r_{E}\prod_{q|N_E, q\neq p}P_q(q^{(p-1)s-1}){\mathcal L}_{E,p}((p-1)s)=\hat{\tau}(\phi(w_{\infty}))\cdot  \hat{\tau}(\iota\phi(\kappa_{\infty})).\]
By {\bf Theorem 2.1} and {\bf Proposition 2.6} the second term is 
\[ \hat{\tau}(\iota\phi(\kappa_{\infty}))=\frac{r_E}{\Omega_E}\prod_{q|N_E, q\neq p}P_q(q^{(p-1)s-1})L^{\dagger}(E,1-(p-1)s).\]
and therefore 
\[{\mathcal L}_{E,p}((p-1)s)=\hat{\tau}(\phi(w_{\infty}))\frac{L^{\dagger}(E,1-(p-1)s)}{\Omega_E}.\]
But what $\hat{\tau}(\phi(w_{\infty}))$ should be? By \cite{Rubin1998} {\bf Lemma A.1 (2)} we know that 
\[ \chi(w_{\infty})=
 \begin{cases}
 W(\chi), & \text{if $\chi$ is nontrivial}\\
 0, & \text{if $\chi$ is trivial.}
 \end{cases}\]
 Now remember that $ W(\chi)$ appears in the functional equation of Dirichlet series:
  \[L(0,\chi)=\frac{1}{\pi i}W(\chi)L(1,\chi^{-1}).\]
 If we were able to replace $\chi$ by $\chi_{s}$  as the introduction it would be
 \[\zeta_{(p)}(s)=\frac{1}{\pi i}\hat{\tau}(\sigma(w_{\infty})) \zeta_{(p)}(1-s),\]
 where $\zeta_{(p)}(s)=(1-p^{-s})\zeta(s)$. In particular
 \[\hat{\tau}(\phi(w_{\infty}))=\frac{\pi i \zeta_{(p)}((p-1)s)}{\zeta_{(p)}(1-(p-1)s)},\]
 and since the zeta function has a simple pole at $s=1$ the order of $\hat{\tau}_{\infty}(\phi(w_{\infty}))$ at $s=0$ is one.
%
\section{Consequences of the conjecture}
\subsection{An application to Mazur-Tate-Teitelbaum conjecture}
Let $E$ be an elliptic curve defined over ${\mathbb Q}$. The following conjecture is due to Mazur, Tate and Teitelbaum.
\begin{conjecture} 
\begin{enumerate}
\item Suppose that $E$ has an good ordinary or a non-split multiplicative reduction at $p$. Then 
 \[{\rm ord}_{s=0}{\mathcal L}_{E,p}(s)={\rm ord}_{s=1}L(E,s).\]
\item Suppose that $E$ has split multiplicative reduction at $p$. Then 
 \[{\rm ord}_{s=0}{\mathcal L}_{E,p}(s)=1+{\rm ord}_{s=1}L(E,s).\]
 \end{enumerate}
\end{conjecture}
Suppose  that $L(E,1)\neq 0$. If $E$ has a good ordinary or a non-split multiplicative reduction at $p$ the conjecture is trivially true by {\bf Fact 2.1}. If $E$ has a split multiplicative reduction it has been proved by Greenberg and Stevens \cite{Greenberg-Stevens}. Later Kobayashi gives an elementary proof  of their statements using Kato's result \cite{Kobayashi2006}. Obviously our conjecture follows from {\bf Conjecture 3.1}. Conversely, using a theorem due to Ono and Skinner \cite{Ono-Skinner1998}, we will show that {\bf Conjecture 1.1} implies {\bf Conjecture 3.1}.\\

Let $\Pi=\{p_1,\cdots, p_t\}$ be a set of mutually distinct primes and take $\epsilon=(\epsilon_1,\cdots \epsilon_t)$ where $\epsilon_{i}\in \{\pm 1\}$. Then we define $P(\Pi,\epsilon)$ to be a set of  square free fundamental discriminants $D$ which satisfy
\[(\frac{D}{p_{i}})=\epsilon_{i}, \quad \forall i,\]
where $(\frac{\cdot}{\cdot})$ denotes Legendre symbol. For a square free fundamental discriminant $D$ let $E_D$ be the twist of $E$ over ${\mathbb Q}(\sqrt{D})$. Namely if 
\[E\,:\, y^2=x^3-ax+b,\quad a,b\in {\mathbb Q},\]
is a Weierstrauss form of $E$, $E_{D}$ is defined to be
\[ E_D\,:\, Dy^2=x^3-ax+b.\]
\begin{fact}(\cite{Ono-Skinner1998} {\bf Corollary 3}) For a positive number $X$ 
\[\sharp\{D\in P(\Pi,\epsilon)\,:\, |D|<X,\, L(E_{D},1)\neq 0\} >> \frac{X}{\log X}.\]
.
\end{fact}
{\bf Proof of {\bf Conjecture 1.1} $\Rightarrow$ {\bf Conjecture 3.1}.} By {\bf Fact 3.1} we know that there is a square free fundamental discriminant $D$ satisfying $(\frac{D}{p})=1$ and $L(E_{D},1)\neq 0$. The first condition guarantees that $E$ and $E_{D}$ have ordinary reduction of the same type at $p$. As we have mentioned before the conjecture is true for $E_{D}$ and {\bf Conjecture 3.1} is derived from {\bf Conjecture 1.1}.
\subsection{A review of Iwasawa theory for an elliptic curve}
Let $E$ be an elliptic curve defined over ${\mathbb Q}$ with good ordinary reduction at a prime $p\geq 5$. 
We denote the Selmer group of $E$ over ${\mathbb Q}_{n}$ by ${\rm Sel}({\mathbb Q}_{n},E[p^{\infty}])$ and define
\[{\rm Sel}({\mathbb Q}_{\infty},E[p^{\infty}]):=\lim_{\rightarrow}{\rm Sel}({\mathbb Q}_{n},E[p^{\infty}]),\]
where limits with respect to restrictions. In general let $M$ be a profinite ${\mathbb Z}_{p}$-module. Its $p$-adic and rational $p$-adic Pontryagin dual is define to be
\[{\mathbb D}(M):={\rm Hom}_{conti}(M, {\mathbb Q}_{p}/{\mathbb Z}_{p}),\]
and
\[{\mathbb D}^{0}(M):={\mathbb D}(M)\otimes_{{\mathbb Z}_p}{{\mathbb Q}_p},\]
respectively. The subscript {\it conti} means the set of continuous homomorphisms. By ${\mathbb Z}_{p}[[{\rm Gal}({\mathbb Q}_{\infty}/{\mathbb Q})]]\simeq {\mathbb Z}_{p}[[t]]$,  we regard $X_{\infty}:={\mathbb D}({\rm Sel}({\mathbb Q}_{\infty},E[p^{\infty}]))$ as a ${\mathbb Z}_{p}[[t]]$-module. It is known that $X_{\infty}$ is a torsion ${\mathbb Z}_{p}[[t]]$-module and therefore ${\mathbb D}^{0}({\rm Sel}({\mathbb Q}_{\infty},E[p^{\infty}]))$ is a torsion $\Lambda_{{\mathbb Q}_{p}}$-module.  Here we put $\Lambda_{{\mathbb Q}_{p}}:={\mathbb Z}_{p}[[t]]\otimes_{{\mathbb Z}_p}{{\mathbb Q}_p}$, which is a discrete valuation ring whose valuation ideal is  generated by $t$. Let ${\rm Char}(X_{\infty})\subset \Lambda_{{\mathbb Q}_{p}}$ be the characteristic ideal of ${\mathbb D}^{0}({\rm Sel}({\mathbb Q}_{\infty},E[p^{\infty}]))$ and the order of its generator with respect to $t$ is denoted by ${\rm ord}_{t}{\rm Char}(X_{\infty})$. 
\begin{prop}
\[{\mathbb D}^{0}({\rm Sel}({\mathbb Q}_{\infty},E[p^{\infty}]))^{{\rm Gal}({\mathbb Q}_{\infty}/{\mathbb Q})}\simeq {\mathbb D}^{0}({\rm Sel}({\mathbb Q},E[p^{\infty}])).\]
\end{prop}
{\bf Proof.} Immediate from \cite{Perrin-Riou1990} {\bf Lemme 6.6}.$\Box$\\

Let $\sf(E/{\mathbb Q})$ be the Shafarevich-Tate group. Taking rational $p$-adic Pontryagin dual of
\[0 \to E({\mathbb Q})\otimes {\mathbb Q}_{p}/{\mathbb Z}_{p} \to {\rm Sel}({\mathbb Q},E[p^{\infty}]) \to \sf(E/{\mathbb Q})[p^{\infty}] \to 0,\]
we have
\[0 \to {\mathbb D}^{0}(\sf(E/{\mathbb Q})[p^{\infty}]) \to {\mathbb D}^{0}({\rm Sel}({\mathbb Q},E[p^{\infty}])) \to {\rm Hom}_{{\mathbb Z}}(E({\mathbb Q}),{\mathbb Z})\otimes{\mathbb Q}_p \to 0.\]
{\bf Propositon 3.1} implies the following theorem.
\begin{thm}
\[{\rm ord}_{t}{\rm Char}(X_{\infty}) \geq {\rm rank}E({\mathbb Q}).\]
Moreover if $\sf(E/{\mathbb Q})[p^{\infty}]$ is finite the equality holds.
\end{thm}

\subsection{Birch and Swinnerton-Dyer conjecture for a semistable elliptic curve}
\begin{lm} Let $E$ be an elliptic curve defined over ${\mathbb Q}$ and $L$ be a quadratic extension of ${\mathbb Q}$. Then for any odd prime $p$ and a positive integer $r$ the restriction gives an isomorphism,
\[\sf(E/{\mathbb Q})[p^r]\simeq \sf(E/L)[p^r].\]
\end{lm}
{\bf Proof.} Since $p$ is odd and since the order of ${\rm Gal}(L/{\mathbb Q})$ is two, $H^{1}({\rm Gal}(L/{\mathbb Q}), E(L))[p^r]=H^{2}({\rm Gal}(L/{\mathbb Q}), E(L))[p^r]=0$. Therefore the inflation-restriction sequence implies
\[H^{1}({\mathbb Q},E)[p^r]\simeq H^{1}(L,E)[p^r].\]
Let $v$ be a place (including $\infty$) of ${\mathbb Q}$. Suppose that $v$ ramifies or inerts in $L$ and let $w$ be the place of $L$ over $v$. The same argument as above shows 
\[H^{1}({\mathbb Q}_{v},E)[p^r]\simeq H^{1}(L_{w},E)[p^r].\]
Suppose $v$ splits in $L$ and let $w$ and $w^{\prime}$ be places of $L$ over $v$. Then both $L_{w}$ and $L_{w^{\prime}}$ are isomorphic to ${\mathbb Q}_{v}$ and 
\[H^{1}({\mathbb Q}_{v},E)[p^r]\stackrel{{\rm Res}}\to H^{1}(L_{w},E)[p^r]\times H^{1}(L_{w^{\prime}},E)[p^r],\]
is the diagonal imbedding,
\[H^{1}({\mathbb Q}_{v},E)[p^r]\stackrel{{\Delta}}\to H^{1}({\mathbb Q}_{v},E)[p^r] \times H^{1}({\mathbb Q}_{v},E)[p^r].\]
 Thus we find that $g$ (reps. $h$) of
 \[\begin{CD}
 0 @>>> \sf(E/{\mathbb Q})[p^r] @>>> H^{1}({\mathbb Q},E)[p^r] @>>> \prod_{v}H^{1}({\mathbb Q}_{v},E)[p^r]\\
 @VVV  @V\text{$f$}VV @V\text{$g$}VV @V\text{$h$}VV\\
 0 @>>>  \sf(E/L)[p^r] @>>>  H^{1}(L,E)[p^r] @>>> \prod_{w}H^{1}(L_{w},E)[p^r],
 \end{CD}\]
is isomorphic (resp. injective). A simple diagram chasing shows that $f$ is an isomorphism.
$\Box$\\
\begin{prop} Let $E$ be an semi-stable elliptic curve defined over ${\mathbb Q}$. Then there is a pair $(D,p)$ such that
\begin{enumerate}
\item $D$ is a square free fundamental discriminant so that
\[L(E_{D},1)\neq 0.\]
\item $p$ is a good ordinary prime of $E$ which satisfies
\[(\frac{D}{p})=1,\quad \sf(E_{D}/{\mathbb Q})[p]=0,\quad p\geq 11.\]
\end{enumerate}
\end{prop}
{\bf Proof.} By {\bf Fact 3.1} there is a fundamental discriminant satisfying (1). Let us fix one of them. Then by \cite{Kato} we know that $\sf(E_{D}/{\mathbb Q})$ is a finite abelian group. Since $E$ is semistable it does not have complex multiplication and, due to Serre, the density of supersingular primes of $E$ is $0$. Therefore there is a prime satisfying (2).
$\Box$\\
The following theorem is a direct consequence of \cite{Skinner-Urban} {\bf Corollary 3.6.10}.
\begin{thm}
Let $E$ be a semistable elliptic curve defined over ${\mathbb Q}$. Let $p\geq 11$ be a prime where $E$ has good ordinary reduction. Then
\[{\rm ord}_{t}{\rm Char}(X_{\infty})={\rm ord}_{s=0}{\mathcal L}_{E,p}(s).\]
\end{thm}
\begin{thm}Let $E$ be a semistable elliptic curve defined over ${\mathbb Q}$. Then we have the following consequences.
\begin{enumerate}
\item The rank of $E({\mathbb Q})$ is equal to ${\rm ord}_{s=1}L(E,s)$.
\item For an odd prime $q$ the $q$-primary part $\sf(E/{\mathbb Q})[q^{\infty}]$ of the Shafarevich-Tate group of $E$ over ${\mathbb Q}$ is finite. Moreover  it is trivial except finitely many primes. 
\end{enumerate}
\end{thm}
{\bf Proof.} Let $(D,p)$ be a pair of {\bf Proposition 3.2}. Since $p$ is a good ordinary prime of $E$ and since $(\frac{D}{p})=1$, $E$ and $E_{D}$ have a good ordinary reduction of the same type at $p$. By $L(E_{D},1)\neq 0$, {\bf Fact 2.1} shows that the order of $p$-adic L-function of $E_{D}$ at the origin is zero. Then by {\bf Conjecture 1.1} and {\bf Theorem 3.2},
\begin{equation}{\rm ord}_{s=1}L(E,s)={\rm ord}_{t}{\rm Char}(X_{\infty}).\end{equation}
On the other hand since $E$ and $E_{D}$ are isomorphic over ${\mathbb Q}(\sqrt{D})$, by {\bf Lemma 3.1},
\begin{equation}\sf(E/{\mathbb Q})[q^r]\simeq \sf(E/{\mathbb Q}(\sqrt{D}))[q^r] \simeq \sf(E_{D}/{\mathbb Q})[q^r],\end{equation}
for any odd prime $q$ and a positive integer $r$. In particular since $\sf(E_{D}/{\mathbb Q})[p]=0$ we see $\sf(E/{\mathbb Q})[p]=0$. Now {\bf Theorem 3.1}  and (14) implies 
\[ {\rm ord}_{s=1}L(E,s)={\rm rank}E({\mathbb Q}).\]
Kato has shown that $L(E_{D},1)\neq 0$ implies finiteness of $ \sf(E_{D}/{\mathbb Q})$(\cite{Kato}). Thus (15) shows that $\sf(E/{\mathbb Q})[q^{\infty}]$ is finite for an odd prime $q$ and moreover vanishes except finitely many primes.
$\Box$
\subsection{Birch and Swinnerton-Dyer conjecture for an elliptic curve defined with a complex multiplication}
Let $E$ be an elliptic curve defined over ${\mathbb Q}$ whose endomorphism ring is isomorphic to the integer ring ${\mathcal O}_{K}$ of a quadratic imaginary field $K$. Then the following is a direct consequence of \cite{Rubin1991}{\bf Theorem 12.3}.

\begin{thm} Let $p\geq 5$ be a prime where $E$ has good ordinary reduction. Then 
\[{\rm ord}_{t}{\rm Char}(X_{\infty})={\rm ord}_{s=0}{\mathcal L}_{E,p}(s).\]
\end{thm}
\begin{thm}Let $E$ be an elliptic curve defined over ${\mathbb Q}$ whose endomorphism ring is isomorphic to ${\mathcal O}_{K}$. Then we have the following consequences.

\begin{enumerate}
\item The rank of $E({\mathbb Q})$ is equal to ${\rm ord}_{s=1}L(E,s)$.
\item For an odd prime $q$ the $q$-primary part $\sf(E/{\mathbb Q})[q^{\infty}]$ of the Shafarevich-Tate group of $E$ over ${\mathbb Q}$ is finite. Moreover it is trivial except finitely many primes. 
\end{enumerate}
\end{thm} 
{\bf Proof.} Let $d_K$ be the discriminant of $K$ and $d_K=q_{1}^{e_1}\cdots q_{t}^{e_t}$ be its factorization by primes. We choose an odd prime $l$ such that $(\frac{d_{K}}{l})=1$. By {\bf Fact 3.1} there is a square free fundamental discriminant $D$ satisfying
\[(\frac{D}{q_{1}})=\cdots (\frac{D}{q_{t}})=(\frac{D}{l})=1,\]
and $L(E_{D},1)\neq 0$. Note that the condition of $\{q_{i}\}_{i}$ implies that $d_K$ and $D$ are coprime. Let $R$ be a multiplicative closed set consisting of $\{-1,0,1\}\subset {\mathbb Z}$  and consider a multiplicative map:
\[{\mathbb Z}/(d_{K}) \times {\mathbb Z}/(D) \stackrel{\chi}\to R\times R, \quad \chi(x,y)=((\frac{d_{K}}{x}),(\frac{D}{y})),\] 
which may be regarded as a map
from ${\mathbb Z}/(d_{K}D)$. Since $\chi(l)=(1,1)$ there are infinitely many primes $q$ satisfying $\chi(q)=(1,1)$ ( In fact it is sufficient that $q\equiv l$ (mod $d_{K}D$)). On the other hand since $L(E_{D},1)\neq 0$, $\sf(E_{D}/{\mathbb Q})$ is a finite group \cite{Rubin1991}. Therefore there is an odd prime $p$ which does not divide the discriminant of $E$ and satisfies
\[\sf(E_{D}/{\mathbb Q})[p]=0, \quad (\frac{d_{K}}{p})=(\frac{D}{p})=1.\]
The last condition gurantees that $E$ has good ordinary reduction at $p$ and that $E$ and $E_{D}$ are of same type at $p$. Now the remaining of a proof is the same as one of {\bf Theorem 3.3} (One uses {\bf Theorem 3.4} in stead of {\bf Theorem 3.2}). $\Box$
\vspace{10mm}


\end{document}